\documentclass[12pt]{article}

\usepackage{amssymb,amsmath,latexsym,amsthm}

\advance \topmargin by -\headheight
\advance \topmargin by -\headsep
\evensidemargin \oddsidemargin
\author{J.-P. Allouche\thanks{The author was partially supported by the ANR
project ``FAN'' (Fractals et Num\'eration), ANR-12-IS01-0002.} \\
CNRS, Institut de Math\'ematiques de Jussieu-PRG \\
\'Equipe Combinatoire et Optimisation \\
Universit\'e Pierre et Marie Curie, Case 247 \\
4 Place Jussieu \\
F-75252 Paris Cedex 05 France \\
{\tt jean-paul.allouche@imj-prg.fr}
}

\title{Paperfolding infinite products and the gamma function}

\date{ }

\def \proof{\bigbreak\noindent{\it Proof.\ \ }}

\def \endpf{{\ \ $\Box$ \medbreak}}

\newtheorem{theorem}{Theorem}
\newtheorem{lemma}{Lemma}
\newtheorem{corollary}{Corollary}
\newtheorem{proposition}{Proposition}
\newtheorem*{conjecture}{Conjecture}

\theoremstyle{definition}

\newtheorem{remark}{Remark}
\newtheorem*{example}{Example}
\newtheorem*{examples}{Examples}

\begin{document}

\maketitle

\begin{abstract}
Taking the product of $(2n+1)/(2n+2)$ raised to the power $+1$ or $-1$ according to the 
$n$-th term of the Thue-Morse sequence gives rise to an infinite product $P$ while replacing 
$(2n+1)/(2n+2)$ with $(2n)/(2n+1)$ yields an infinite product $Q$, where
$$
P = \left(\frac{1}{2}\right)^{+1} \left(\frac{3}{4}\right)^{-1} \left(\frac{5}{6}\right)^{-1}
\left(\frac{7}{8}\right)^{+1}\cdots
\ \mbox{\rm and} \ \
Q = \left(\frac{2}{3}\right)^{+1} \left(\frac{4}{5}\right)^{-1} \left(\frac{6}{7}\right)^{-1}
\left(\frac{8}{9}\right)^{+1}\cdots
$$
%with $P = (1/2)(4/3)(6/5)(7/8)(10/9)\cdots$ and $Q = (2/3)(5/4)(7/6)(8/9)(11/10)\cdots$. 
Though it is known that $P = 2^{-1/2}$, nothing is known about $Q$. Looking at the corresponding 
question when the Thue-Morse sequence is replaced by the regular paperfolding sequence, we obtain two 
infinite products $A$ and $B$, where
$$
A = \left(\frac{1}{2}\right)^{+1} \left(\frac{3}{4}\right)^{+1} \left(\frac{5}{6}\right)^{-1}
\left(\frac{7}{8}\right)^{+1}\cdots
\ \mbox{\rm and} \ \
B = \left(\frac{2}{3}\right)^{+1} \left(\frac{4}{5}\right)^{+1} \left(\frac{6}{7}\right)^{-1}
\left(\frac{8}{9}\right)^{+1}\cdots
$$
%$A = (1/2)(3/4)(6/5)(7/8)(9/10)\cdots$ and $B = (2/3)(4/5)(7/6)(8/9)(10/11)\cdots$.
Here nothing is known for $A$, but we give a closed form for $B$ that involves the value of 
the gamma function at $1/4$. We then prove general results where $(2n+1)/(2n+2)$ or $(2n)/(2n+1)$
are replaced by specific rational functions. The corresponding infinite products have a closed 
form involving gamma values. In some cases there is no explicit gamma value occurring in
the closed-form formula, but only trigonometric functions.

\medskip


\medskip

\noindent
{\bf MSC Classes}: 11B85, 11A63, 11J81, 33B15, 68R15.

\end{abstract}

\section{Introduction}

Unexpected explicit values for infinite series or infinite products are somehow
fascinating. One of the most famous examples is the celebrated {\it tour de force} 
of Euler addressing the Basel (or Mengoli) problem and finding the closed form 
$$
\frac{1}{1^2} + \frac{1}{2^2} + \frac{1}{3^2} + \cdots = \frac{\pi^2}{6}\cdot
$$
This formula is indeed unexpected: how could we guess that the square of the area of a 
circle with radius $1$ would be occurring here? Of course this suggests looking at the more 
general sums $\zeta(k) = \sum_n \frac{1}{n^k}$: we know that Euler himself gave their values 
(also unexpected) for $k$ an even integer, while the case where $k$ is odd is still largely open, 
even after Ap\'ery's proof of the irrationality of $\zeta(3)$ (see \cite{vdp}), Ball-Rivoal's proof 
that infinitely many $\zeta(k)$ with $k$ odd are irrational \cite{BR, Rivoal}, and Zudilin's proof 
that at least one of the reals $\zeta(5), \zeta(7), \zeta(9), \zeta(11)$ is irrational \cite{Zudilin}. 

\bigskip

Another classical closed formula, which also looks intriguing at first sight, is the infinite product
discovered by Wallis about 80 years before Euler's result, namely
$$
\frac{\pi}{2} = \frac{2}{1} \times \frac{2}{3} \times \frac{4}{3} \times \frac{4}{5}
\times \frac{6}{5} \times \frac{6}{7} \times \frac{8}{7} \times \frac{8}{9} \times \cdots 
= \prod_{n=1}^{\infty} \frac{(2n)(2n)}{(2n-1)(2n+1)}\cdot
$$
Again there is no immediate intuition why $\pi$ should occur when multiplying the even 
numbers twice and dividing by the odd numbers twice.

\bigskip

Much less famous, but unexpected as well, is the Woods-Robbins infinite product \cite{Woods, Robbins}. 
Define $s(n)$ to be the sum of the binary digits of the integer $n$, and consider the sequence
$(-1)^{s(n)}$. In other words $m_n = (-1)^{s(n)}$ can be defined by
$$
m_0 = 1 \ \mbox{\rm and, for all $n \geq 0$,} \
m_{2n} = m_n \ \mbox{\rm and} \ m_{2n+1} = -m_n.
$$
This sequence is known as the (Prouhet-)Thue-Morse sequence, and its first terms are 
$+1, -1, -1, +1, \ldots$ Then
\begin{equation}\label{TM}
\left(\frac{1}{2}\right)^{+1} \left(\frac{3}{4}\right)^{-1} \left(\frac{5}{6}\right)^{-1} 
\left(\frac{7}{8}\right)^{+1} \cdots =
\prod_{n \geq 0} \left(\frac{2n+1}{2n+2}\right)^{(-1)^{s(n)}} = \frac{\sqrt{2}}{2}\cdot
\end{equation}

More general products of the form $\prod_n R(n)^{u_n}$ where $(u_n)_{n \geq 0}$ is a sequence with 
``regularity properties'', and $R(n)$ is an {\it ad hoc} rational function can be looked at.
Among possible sequences $(u_n)_{n \geq 0}$, one can think of {\it automatic sequences}:
a sequence $(u_n)_{n \geq 0}$ is $q$-automatic for some integer $q \geq 2$ if its $q$-kernel, i.e.,
the set of subsequences $\{(u_{q^k n + j})_{n \geq 0}, \ k \geq 0, \ j \in [0, q^k-1]\}$ is finite
(see, e.g., \cite{AS} to read more on automatic sequences). While finding closed forms of 
$\prod_n R(n)^{u_n}$ where $(u_n)_{n \geq 0}$ is any automatic sequence seems out of reach, closed 
forms for particular $2$-automatic sequences were found (see, e.g., \cite{AllCohen, ASlms, AllSon, LouPro}). 
These are typically sequences $((-1)^{w_n})_{n \geq 0}$ where $w_n$ counts the number of occurrences 
of some pattern in the binary expansion of $n$: for example the Thue-Morse sequence 
$((-1)^{s(n)})_{n \geq 0}$ (the sum of the binary digits of $n$ is also the number of $1$'s in the binary
expansion of $n$). Another classical, but totally different, $2$-automatic sequence is the regular 
$\pm 1$ paperfolding sequence $\varepsilon_n$ which is obtained by iteratively folding a rectangular 
piece of paper, see, e.g., \cite[p.\ 155--156]{AS}. The regular paperfolding sequence can be defined by
$$
\mbox{\rm for all $n \geq 0$,} \
\varepsilon_{2n} = (-1)^n \ \mbox{\rm and} \ \varepsilon_{2n+1} = \varepsilon_n.
$$

If we try to mimic one proof (see \cite[Proposition~5, p.\ 6]{ubiq}) of Equality~(\ref{TM}) above,
we can proceed as follows. Let $A$ and $B$ be the two infinite products (they are convergent
as will be seen in Proposition~\ref{conv} below) defined by
$$
A := \prod_{n \geq 0} \left(\frac{2n+1}{2n+2}\right)^{\varepsilon_n} 
\ \ \ \ \ \mbox{\rm and} \ \ \ \ \
B := \prod_{n \geq 1} \left(\frac{2n}{2n+1}\right)^{\varepsilon_n}\cdot
$$
Multiplying these two products, we can write
$$
AB = \frac{1}{2} \prod_{n \geq 1} \left(\frac{(2n+1)(2n)}{(2n+2)(2n+1)}\right)^{\varepsilon_n}
= \frac{1}{2} \prod_{n \geq 1} \left(\frac{n}{n+1}\right)^{\varepsilon_n}\cdot
$$
Splitting the product on the right into odd and even indexes, we thus have
$$
AB = \frac{1}{2} \prod_{n \geq 1} \left(\frac{2n}{2n+1}\right)^{\varepsilon_{2n}}
\prod_{n \geq 0} \left(\frac{2n+1}{2n+2}\right)^{\varepsilon_{2n+1}} =
\frac{1}{2} \prod_{n \geq 1} \left(\frac{2n}{2n+1}\right)^{(-1)^n}
\prod_{n \geq 0} \left(\frac{2n+1}{2n+2}\right)^{\varepsilon_n}\cdot
$$
Since the last product is $A$ ($\neq 0$), we have 
$$
B = \frac{1}{2} \prod_{n \geq 1} \left(\frac{2n}{2n+1}\right)^{(-1)^n} 
= \frac{1}{2} \prod_{n \geq 1} \left(\frac{4n}{4n+1}\right) 
\prod_{n \geq 0} \left(\frac{4n+3}{4n+2}\right) =
\frac{1}{2} \prod_{n \geq 0} \left(\frac{(4n+4)(4n+3)}{(4n+5)(4n+2)}\right)\cdot
$$
Using classical results on the gamma function (Proposition~\ref{ww} below and the
reflection formula) we obtain
\begin{equation}\label{PP}
\prod_{n \geq 1} \left(\frac{2n}{2n+1}\right)^{\varepsilon_n} =
\frac{\Gamma(1/4)^2}{8\sqrt{2 \pi}}\cdot
\end{equation} 

The purpose of this paper is to show how to obtain closed forms for certain infinite 
products $\prod_n R(n)^{\varepsilon_n}$, where $R(n)$ belongs to some families of 
rational functions. These closed forms involve, as in Equality~(\ref{PP}), values of 
the $\Gamma$ function (and of trigonometric functions).

\section{Two premiminary results}

In this section we give three preliminary results. The first two are classical.

\begin{proposition}{\rm (see, e.g., \cite[Section~12-13, p.\ 238--239]{WW})}\label{ww}
Let $d$ be a positive integer. Let $(a_i)_{1 \leq i \leq d}$ and $(b_j)_{1 \leq j \leq d}$
be complex numbers such that no $a_i$ and no $b_j$ belongs to $\{0, -1, -2, \ldots\}$.
If $a_1 + a_2 + \cdots + a_d = b_1 + b_2 + \cdots + b_d$, then
$$
\prod_{n \geq 0}\frac{(n+a_1) \cdots (n+a_d)}{(n+b_1) \cdots (n+b_d)}
= \frac{\Gamma(b_1) \cdots \Gamma(b_d)}{\Gamma(a_1) \cdots \Gamma(a_d)}\cdot
$$
\end{proposition}

\begin{proposition}{\rm (see, e.g., \cite[Sections~12-14 and 12-15, p.\ 239--240]{WW})}\label{gamma-rel}
The Gamma function satisfies the ``reflection formula'' and the ``duplication formula'' respectively given 
by
$$
\Gamma(z)\Gamma(1-z)= \pi/\sin(\pi z) \ \ \mbox{\rm and} \ \ 
2^{2z-1} \Gamma(z) \Gamma(z + \frac{1}{2}) = \sqrt{\pi} \Gamma(2z).
$$

\end{proposition}

The third result that we need is the asymptotic behavior of the summatory function of
the (regular) $\pm 1$ paperfolding sequence and an application to the convergence of certain series.

\begin{proposition}\label{conv}
Let $(\varepsilon_n)_{n \geq 0}$ be the $\pm 1$ paperfolding sequence. 

\begin{itemize}

\item[(i)]
We have the upper bound
$$
\sum_{0 \leq k < n} \varepsilon_k = O(\log n).
$$

\item[(ii)] Let $f$ be a map from ${\mathbb R}$ to ${\mathbb C}$ such that, when $x$ tends to infinity,
$f(x)$ tends to zero and $|f(x+1) - f(x)| = O(1/x^a)$ for some $ a > 1$. Then the series 
$\sum \varepsilon_n f(n)$ is convergent.

\end{itemize}
\end{proposition}

\proof 

(i) Let $S(n) = \sum_{0 \leq k < n} \varepsilon_k$. Looking at $S(2n)$ and splitting the sum
into even and odd indexes yields
\begin{equation}\label{S2n}
\begin{array}{lll}
S(2n) = \displaystyle\sum_{0 \leq k < 2n} \varepsilon_k &=&
\displaystyle\sum_{0 \leq k < n} \varepsilon_{2k} + \sum_{0 \leq k < n} \varepsilon_{2k+1} =
\sum_{0 \leq k < n} (-1)^k + \sum_{0 \leq k < n} \varepsilon_k \\
&=& \displaystyle\sum_{0 \leq k < n} (-1)^k + S_n = \frac{1- (-1)^n}{2} + S_n.
\end{array}
\end{equation}

Let us prove that for all $n \geq 1$ we have $|S(n)| \leq 1 + \log_2(n)$ (where $\log_2$ is the base 
$2$ logarithm). Since $S(1) = \varepsilon_0 = 1 = 1+ \log_2 1$, it suffices to prove that if the inequality 
is true for $n$, then it is true for both $2n$ and $2n+1$. (Hint: this implies by induction on $N$ that, 
if the property is true on $[1, 2^N-1]$ then it is true on $[1, 2^{N+1}-1]$.) Equality~(\ref{S2n})
above implies
$$
\begin{array}{lll}
|S_{2n}| &\leq& 1 + |S_n| \leq 2 + \log_2(n) = 1 + \log_2(2n) \\
|S_{2n+1}| &=& |S_{2n} + (-1)^n| = \displaystyle\frac{1+(-1)^n}{2} + S_n \\
&\leq& 1 + |S_n| \leq 2 + \log_2(n) = \log_2(2n) < \log_2(2n+1).\\
\end{array}
$$

\medskip

(ii) This is an immediate consequence of (i) by Abel summation. \endpf

\section{A general product involving the paperfolding sequence}

We give here a general result, involving the $\pm 1$ paperfolding sequence. (We use the
notation ${\mathbb R}^+$ for the set of non-negative real numbers and ${\mathbb R}^{+*}$ 
for the set of positive real numbers.)

\begin{lemma}\label{general}
Let $g$ be a map from ${\mathbb R}^+$ to ${\mathbb R}^{+*}$ such that, when $x$ tends to infinity, 
$g(x)$ tends to $1$ and $g(x+1)/g(x) = 1 + O(1/x^a)$ for some $a > 1$. Then 
$$
\prod_{n \geq 0} \left(\frac{g(n)}{g(2n+1)}\right)^{\varepsilon_n} 
= \prod_{n \geq 0} \frac{g(4n)}{g(4n+2)}\cdot
$$
\end{lemma}

\proof Convergence of the infinite products results from Proposition~\ref{conv} (ii) with $f = \log g$.
We then write
$$
\prod_{n \geq 0} g(n)^{\varepsilon_n} = \prod_{n \geq 0} g(2n)^{\varepsilon_{2n}}
                                        \prod_{n \geq 0} g(2n+1)^{\varepsilon_{2n+1}}
= \prod_{n \geq 0} g(2n)^{(-1)^n} \prod_{n \geq 0} g(2n+1)^{\varepsilon_{2n+1}}.
$$
Hence
$$
\prod_{n \geq 0} \left(\frac{g(n)}{g(2n+1)}\right)^{\varepsilon_n} = \prod_{n \geq 0} g(2n)^{(-1)^n}
= \prod_{n \geq 0} \frac{g(4n)}{g(4n+2)}\cdot \ \Box
$$

\begin{example}
Taking $g : {\mathbb R}^+ \to {\mathbb R}^{+*}$ defined by $g(x) = \frac{n}{n+1}$ if $x > 0$ 
and $g(0) = 1$, and using Propositions~\ref{ww} and \ref{gamma-rel}, we obtain 
Equality~(\ref{PP}) of the introduction:
$$
\prod_{n \geq 1} \left(\frac{2n}{2n+1}\right)^{\varepsilon_n} = \frac{\Gamma(1/4)^2}{8\sqrt{2 \pi}}\cdot
$$
\end{example}

\begin{remark}
Using a classical result on elliptic integrals (see, e.g., \cite[p.\ 373]{ChSe}),
it is amusing though anecdotical to note that
$$
\prod_{n \geq 1} \left(\frac{2n}{2n+1}\right)^{\varepsilon_n} = 
\frac{1}{2} \int_0^{\pi/2} \frac{\mbox{\rm d}\varphi}{\sqrt{2 - \sin^2 \varphi}}\cdot
$$

\end{remark}

\section{Paperfolding infinite products and the gamma function}

In this section we give applications of Lemma~\ref{general} above to finding closed-form 
expressions for ``simple'' infinite products involving the paperfolding sequence. 

\begin{theorem}\label{gamma-simple}
Let $b$ be a positive real number. We have the following expression.
$$
\prod_{n \geq 0} \left( \frac{n+b}{n + \frac{1+b}{2}} \right)^{\varepsilon_n}
= \frac{\Gamma(\frac{1}{4})^2}{\pi \sqrt{2}} \times \frac{\Gamma(\frac{1}{2} + \frac{b}{4})}
       {\Gamma(\frac{b}{4})} = 
2^{\frac{1-b}{2}} \left(\frac{\Gamma(\frac{1}{4})}{\Gamma(\frac{b}{4})}\right)^2 
\times \frac{\Gamma(\frac{b}{2})}{\sqrt{\pi}}\cdot
$$
\end{theorem}

\proof Apply Lemma~\ref{general} with $g(x) = \frac{x+b}{x+1}$ and Propositions~\ref{ww} 
and \ref{gamma-rel}.  \endpf

% VOIR http://arxiv.org/pdf/1212.0251v5.pdf

\begin{examples}

Taking $b=2$ (resp. $b=3$) in Theorem~\ref{gamma-simple} above, we obtain
$$
\prod_{n \geq 0} \left(\frac{2n+4}{2n+3}\right)^{\varepsilon_n} = \frac{\Gamma(1/4)^2}{\pi^{3/2} \sqrt{2}}
\ \ \mbox{\rm and} \ \
\prod_{n \geq 0} \left(\frac{n+3}{n+2}\right)^{\varepsilon_n} = \frac{\Gamma(1/4)^4}{8\pi^2}\cdot
$$

\end{examples}

\bigskip

\begin{remark} A totally unexpected relation can be obtained from Theorem~\ref{gamma-simple} above
with $b=1+\frac{2x}{\pi}$ and from Relation~(14) in the paper \cite{Blagouchine} of Blagouchine
$$
\prod_{n \geq 0} \left(\frac{n + 1 + \frac{2x}{\pi}}{n + 1 + \frac{x}{\pi}} \right)^{\varepsilon_n}
= \frac{1}{2\pi} \int_0^{\infty} \frac{\log(1+\frac{x^2}{u^2})}{\cosh u} \, \mbox{\rm d}u.
$$
\end{remark}

\begin{remark}\label{k-ell}
Taking a nonnegative integer $k$, and replacing $b$ by $\frac{b + 2^k - 1}{2^k}$ in 
Theorem~\ref{gamma-simple} above, we get
$$
T_k(b) = \prod_{n \geq 0} \left(\frac{n + \frac{b + 2^k - 1}{2^k}}{n + \frac{b + 2^{k+1} -1}{2^{k+1}}}\right)^{\varepsilon_n}
= \frac{\Gamma(\frac{1}{4})^2}{\pi \sqrt{2}} \times \frac{\Gamma(\frac{3 \times 2^k - 1 + b}{2^{k+2}})}
       {\Gamma(\frac{2^k - 1 + b}{2^{k+2}})} =
2^{\frac{1-b}{2^{k+1}}} \left(\frac{\Gamma(\frac{1}{4})}{\Gamma(\frac{2^k-1+b}{2^{k+2}})}\right)^2
\times \frac{\Gamma(\frac{2^k - 1 + b}{2^{k+1}})}{\sqrt{\pi}}\cdot
$$
By multiplying together these formulas for consecutive values of $k$, we can obtain a closed 
formula for the infinite products $T_{k,\ell}(b)$ (where $k$ and $\ell$ are nonnegative integers
with $k > \ell$, and where $b$ is a positive real) 
$$
T_{k,\ell}(b) = \prod_{n \geq 0} \left(\frac{n + \frac{b + 2^k - 1}{2^k}}
                        {n + \frac{b + 2^{\ell} -1}{2^{\ell}}}\right)^{\varepsilon_n}
= T_k(b) T_{k+1}(b) \cdots T_{\ell-1}(b).
$$
Thus, by replacing in $T_{k,\ell}(b)$ the quantity $\frac{b + 2^{\ell} -1}{2^{\ell}}$ by $c$, we have
a closed form formula for the infinite product $U_j(c)$, defined by (where $j$ is a nonnegative integer,
and where $c$ is a real $> 1 - \frac{1}{2^j}$)
$$
U_j(c) = \prod_{n \geq 0} \left(\frac{n + 2^j c + 1 - 2^j}{n + c}\right)^{\varepsilon_n}\cdot
$$
\end{remark}

\begin{example}
An example of an infinite paperfolding product obtained with Remark~\ref{k-ell} is
$$
\prod_{n \geq 0} \left(\frac{2n+6}{2n+3}\right)^{\varepsilon_n} = 
\frac{\Gamma(1/4)^6}{8\sqrt{2}\,\pi^{7/2}}\cdot
$$
Namely
$$
\prod_{n \geq 0} \left(\frac{2n+6}{2n+3}\right)^{\varepsilon_n} = 
\prod_{n \geq 0} \left(\frac{n+3}{n+\frac{3}{2}}\right)^{\varepsilon_n} = T_{0,2}(3) = T_0(3) T_1(3)
= \frac{\Gamma(1/4)^4}{8\pi^2} \frac{\Gamma(1/4)^2}{\sqrt{2} \pi^{3/2}}
= \frac{\Gamma(1/4)^6}{8\sqrt{2}\,\pi^{7/2}}\cdot
$$

\end{example}

\begin{remark}
We do not know how to characterize the positive integers $u, v, w$ for which the infinite product 
$\displaystyle\prod_{n \geq 0} \left(\frac{un+v}{un+w}\right)^{\varepsilon_n}$ can be expressed
in a closed form (also see Remark~\ref{do-not-know}).
\end{remark}

\bigskip

We can give a more symmetric expression which will have the advantage of yielding ``gamma-free'' 
values of paperfolding infinite products.

\begin{theorem}\label{gamma}
Let $b$ and $c$ be two positive real numbers. Then
$$
\prod_{n \geq 0} \left(
\frac{(n+b)(n + \frac{1+c}{2})}{(n+c)(n + \frac{1+b}{2})}
                 \right)^{\varepsilon_n} 
= \frac{\Gamma(\frac{c}{4})\Gamma(\frac{1}{2} + \frac{b}{4})}
       {\Gamma(\frac{b}{4})\Gamma(\frac{1}{2} + \frac{c}{4})}\cdot
$$
\end{theorem}

\proof Apply Theorem~\ref{gamma-simple} for $b$ and $c$ and compute the quotient. \endpf

\begin{remark}
Another possibly anecdotal remark is that the right side of the equality in Theorem~\ref{gamma}
above can be seen as a coefficient in a connection formula for algebraic hypergeometric
functions (see, e.g., \cite[p.\ 107]{Erdelyi}; also see \cite{Kato1, Kato2}): 
with the notation of \cite{Kato1},
$\alpha_1^2(u, v, w) := \frac{\Gamma(w) \Gamma(w-u-v)}{\Gamma(w-u)\Gamma(w-v)}$, we have
$$
\prod_{n \geq 0} \left(\frac{(n+b)(n + \frac{1+c}{2})}{(n+c)(n + \frac{1+b}{2})}
                 \right)^{\varepsilon_n} =
\alpha_1^2\left(\frac{1}{2}, \frac{b-c}{4}, \frac{1}{2} + \frac{b}{4}\right)\cdot 
$$

\end{remark}

\section{``Gamma-free'' paperfolding products}

In this section we describe ``gamma-free'' paperfolding infinite products, i.e., products 
where the gamma function ``does not appear explicitly''. We begin with a ``trigonometric''
corollary of Theorem~\ref{gamma} above.

\medskip

\begin{corollary}\label{tangent}

\ { } 

(i) Let $b$ be a real number in $(0, 2)$. Then
$$
\prod_{n \geq 0} \left(\frac{(n+b)(2n+3-b)}{(n+2-b)(2n+1+b)}\right)^{\varepsilon_n}
= \tan \frac{\pi b}{4}\cdot
$$

(ii) Let $k$ be an integer $\geq 3$. Then

$$
\begin{array}{lll}

\displaystyle\prod_{n \geq 0} \left(\frac{(kn+k-1)(2kn+2k+1)}{(kn+k+1)(2kn+2k-1)}\right)^{\varepsilon_n} 
&=& \displaystyle\tan \frac{(k-1)\pi}{4k} \\
\displaystyle\prod_{n \geq 0} \left(\frac{(kn+k-2)(kn+k+1)}{(kn+k+2)(kn+k-1)}\right)^{\varepsilon_n}
&=& \displaystyle\tan \frac{(k-2)\pi}{4k} \\
\displaystyle\prod_{n \geq 0} \left(\frac{(kn+k-2)(2kn+2k+1)}{(kn+k+2)(2kn+2k-1)}\right)^{\varepsilon_n}
&=& \displaystyle\tan \frac{(k-1)\pi}{4k} \tan \frac{(k-2)\pi}{4k}\cdot

\end{array}
$$
\end{corollary}

\proof Relation~(i) is a consequence of Theorem~\ref{gamma} above with $c = 2-b$, and of the
reflection formula $\Gamma(z)\Gamma(1-z)= \pi/\sin(\pi z)$ (Proposition~\ref{gamma-rel} above). 
The first two relations in (ii) are obtained from (i) by taking $b = (k-1)/k$ and $b = (k-2)/k$; 
multiplying them together gives the last relation. \endpf

\bigskip
\medskip

Relations obtained for some particular values of $b$ or of $k$ are of interest. In particular 
it has been known since Lambert (see \cite[pp.~133--139]{Lambert}; or see, e.g., \cite{Garibaldi}) 
that trigonometric functions at angles whose values in degrees are $0, 3, 6, 9, \ldots$ can be 
written using only addition, subtraction, multiplication, division, and (possibly nested) square 
roots, of rational numbers. Using these expressions in (i) and then in (ii) yields, e.g., the 
following paperfolding infinite products.

\bigskip
\bigskip

\begin{examples}

\ { }

\begin{itemize}

\item For $b = \frac{3}{2}$ we obtain
$$
\prod_{n \geq 0} \left(\frac{(2n+3)(4n+3)}{(2n+1)(4n+5)}\right)^{\varepsilon_n}
= \tan \frac{3\pi}{8} = \cot \frac{\pi}{8} = 1 + \sqrt{2}.
$$

\item For $b = \frac{1}{5}$ we obtain
$$
\prod_{n \geq 0} \left(\frac{(5n+1)(5n+7)}{(5n+9)(5n+3)}\right)^{\varepsilon_n}
= \tan \frac{\pi}{20} = \sqrt{5} + 1 - \sqrt{5 + 2\sqrt{5}}.
$$

\item For $b = \frac{2}{5}$ we obtain
$$
\prod_{n \geq 0} \left(\frac{(5n+2)(10n+13)}{(5n+8)(10n+7)}\right)^{\varepsilon_n}
= \tan \frac{\pi}{10} = \frac{1}{5} \sqrt{5(5 - 2 \sqrt{5})}.
$$

\item For $k=3$ we obtain
$$
\begin{array}{lll}
\displaystyle\prod_{n \geq 0} \left(\frac{(3n+2)(6n+7)}{(3n+4)(6n+5)}\right)^{\varepsilon_n}
&=& \displaystyle\tan \frac{\pi}{6} = \frac{\sqrt{3}}{3}
\\
\displaystyle\prod_{n \geq 0} \left(\frac{(3n+1)(3n+4)}{(3n+5)(3n+2)}\right)^{\varepsilon_n}
&=& \displaystyle\tan \frac{\pi}{12} = 2 - \sqrt{3}
\\
\displaystyle\prod_{n \geq 0} \left(\frac{(3n+1)(6n+7)}{(3n+5)(6n+5)}\right)^{\varepsilon_n}
&=& \displaystyle\tan \frac{\pi}{6} \tan \frac{\pi}{12} = \frac{2\sqrt{3}}{3} - 1.
\end{array}
$$

\item For $k=5$ we obtain

$$
\begin{array}{lll}
\displaystyle\prod_{n \geq 0} \left(\frac{(5n+3)(5n+6)}{(5n+7)(5n+4)}\right)^{\varepsilon_n}
&=& \displaystyle\tan \frac{3\pi}{20} = \sqrt{5} - 1 - \sqrt{5 - 2 \sqrt{5}}
\\
\displaystyle\prod_{n \geq 0} \left(\frac{(5n+4)(10n+11)}{(5n+6)(10n+9)}\right)^{\varepsilon_n}
&=& \displaystyle\tan \frac{\pi}{5} = \sqrt{5 - 2 \sqrt{5}}
\\
\displaystyle\prod_{n \geq 0} \left(\frac{(5n+3)(10n+11)}{(5n+7)(10n+9)}\right)^{\varepsilon_n}
&=& \sqrt{25 - 10 \sqrt{5}} - \sqrt{5 - 2 \sqrt{5}} - 5 + 2\sqrt{5}.
\end{array}
$$

\end{itemize}

\end{examples}

\bigskip

One may ask whether it is possible to find other cases of ``gamma-free'' paperfolding products. In 
particular the so-called ``short gamma products'' may give paperfolding products where the gamma 
function does not occur explicitly. Results on short gamma products are due to S\'andor and T\'oth 
\cite{ST} and to Nijenhuis \cite{Nijenhuis}, while ``isolated'' examples can be found in the literature 
(see the paper of Borwein and Zucker \cite{BZ}; also see \cite{Nijenhuis}). Before stating these 
results we give some notation: $\varphi$ is the Euler function; the set $\Phi(m)$ is the set of integers in 
$[0, m]$ that are relatively prime to $m$. We note that $\Phi(m)$ is a group with respect to multiplication 
modulo $m$.

\begin{theorem}\label{short}

\ {  }

\begin{itemize}

\item[(i)] \mbox{\rm (S\'andor and T\'oth \cite{ST}; also see \cite{ChSt, Martin})}
The following relation holds.
$$
\prod_{\stackrel{{\scriptstyle 1 \leq k \leq m}}{\gcd(k,m) = 1}} 
\Gamma\left(\frac{k}{m}\right) = 
\begin{cases}
\displaystyle\frac{(2\pi)^{\varphi(m)/2}}{\sqrt{p}}, \ \mbox{if $m$ is a prime power;} \\
(2\pi)^{\varphi(m)/2}, \ \mbox{otherwise}.
\end{cases}
$$

\item[(ii)] \mbox{\rm (Nijenhuis \cite{Nijenhuis})}
Let $n > 1$ be an odd integer. Let $A_n$ be the cyclic subgroup of $\Phi(2n)$ generated by $(n+2)$ or
any one of its cosets. Let $\nu(n)$ denote the cardinality of $A_n$, and $b(A_n)$ the number of elements
of $A_n$ that are larger than $n$. Then
$$
\prod_{x \in A_n} \Gamma\left(\frac{x}{2n}\right) = 2^{b(A)} \pi^{\nu(n)/2}.
$$

\item[(iii)] \mbox{\rm (Borwein and Zucker \cite{BZ}; also see \cite[section~5]{Nijenhuis})}
The following relations hold.
$$
\begin{array}{lll}
\displaystyle\Gamma\left(\frac{1}{24}\right)\Gamma\left(\frac{11}{24}\right)
\Gamma\left(\frac{19}{24}\right)\Gamma\left(\frac{17}{24}\right)
&=& 4 \pi^2 \, 5^{1/4}\\
&&\\
\displaystyle\Gamma\left(\frac{1}{20}\right)\Gamma\left(\frac{9}{20}\right)
\Gamma\left(\frac{13}{20}\right)\Gamma\left(\frac{17}{20}\right)
&=& 4 \pi^2 \, \sqrt{3}.
\end{array}
$$

\end{itemize}
\end{theorem}

\begin{remark}
The ``sporadic'' formulas in Theorem~\ref{short}~(iii) can be retrieved from the relations for the 
values $\Gamma(p/q)$ with $q | 60$ given in \cite{Vidunas}.
\end{remark}
 
In order to use such a ``short'' gamma product in the infinite products given in Theorem~\ref{gamma}
we would like to have four terms in the theorem above (up to multiplying the denominator
$\Gamma(\frac{b}{4})\Gamma(\frac{1}{2} + \frac{c}{4})$ in Theorem~\ref{gamma} by 
$\Gamma(1-\frac{b}{4})\Gamma(\frac{1}{2} - \frac{c}{4})$ and using the reflection formula).
Since the number of terms in the short product of Theorem~\ref{short}~(i) is $\varphi(n)$, we need
to restrict to the integers $m$ such that $\varphi(m) = 4$ (which is easily proven equivalent
to $m \in \{5, 8, 10, 12\}$). We did not obtain new products this way. Using Theorem~\ref{short}~(ii)
gives new results but they are not gamma-free. Finally we can deduce from Theorem~\ref{short}~(iii)
the following result.

\begin{theorem}\label{sporadic-appl}
The following relations hold.
$$
\begin{array}{lll}
\displaystyle\prod_{n \geq 0} \left(\frac{(6n+5)(12n+7)}{(6n+1)(12n+11)}\right)^{\varepsilon_n}
&=& \displaystyle 4 \times 5^{1/4} \times \sin \frac{5\pi}{24} \times \sin \frac{11\pi}{24} 
= 5^{1/4}(1+\sqrt{2}). \\
\displaystyle\prod_{n \geq 0} \left(\frac{(5n+3)^2}{(5n+1)(5n+4)}\right)^{\varepsilon_n}
&=& \displaystyle 4 \sqrt{3} \times \sin \frac{3\pi}{20} \times \sin \frac{11\pi}{20}
= \frac{\sqrt{6}\sqrt{5-\sqrt{5}} + \sqrt{15} - \sqrt{3}}{2}\cdot
\end{array}
$$
\end{theorem}

\proof Put $b=5/6$ and $c=1/6$ (resp.\ $b=3/5$ and $c = 1/5$) in Corollary~\ref{gamma}.
Use Theorem~\ref{short}~(iii) and the reflection formula. \endpf

\section{Arithmetical nature of the paperfolding products} 

So far trying to give closed-form expressions for infinite products of the type $\prod_n R(n)^{u_n}$ 
where $(u_n)_{n \geq 0}$ is some automatic sequence, is possible when $u_n = (-1)^{w_n}$ (where $w_n$ 
counts the number of occurrences of some pattern in the binary expansion of $n$) as recalled in the
introduction, or, as done above, when $u_n = \varepsilon_n$ is the $\pm 1$ paperfolding sequence.
But the choice for the corresponding $R(n)$'s is in both cases drastically limited. This does not apply
to the case where $u_n = 1$ for all $n$; see Proposition~\ref{ww} where the only limitation is that the 
infinite product converges. 

\bigskip

A partly related question is the arithmetical nature of these infinite products: when are they algebraic, 
and where are they transcendental? The question is certainly not easy, even when we have a closed-form 
expression. In the examples above, we know that: 

\begin{itemize}

\item all infinite products given in Corollary~\ref{tangent}~(i) for $b$ rational, and in 
      Corollary~\ref{tangent}~(ii) are algebraic; 

\item both products in Theorem~\ref{sporadic-appl} are algebraic;

\item the two examples of products following Theorem~\ref{gamma-simple} as well as the example 
      following Remark~\ref{k-ell} are transcendental from a result of \v{C}udnovs'ki\u{\i} 
      (see \cite{Cudnovskii}; also see, e.g., \cite[Theorem~14, p.\ 441]{Waldschmidt}) stating 
      that the numbers $\Gamma(1/4)$ and $\pi$ are algebraically independent;

\end{itemize}

\begin{remark}    

In order to prove the transcendence of the infinite product (as previously $(\varepsilon_n)_n$ is the
paperfolding sequence) $\prod_{n \geq 1} \left(\frac{2n}{2n+1}\right)^{\varepsilon_n} = 
\frac{\Gamma(1/4)^2}{8\sqrt{2 \pi}}$, we do not need the full strength of \v{C}udnovs'ki\u{\i}'s result: 
it suffices to use, e.g., \cite[Corollary 7.4, p.\ 200]{BT}: the transcendence of $\Gamma^2(1/4)/\sqrt{\pi}$ 
can be proved by relating this constant to a nonzero period of the Weierstrass ${\mathcal P}$-function, i.e., 
of the elliptic curve $y^2 = 4x^3 - 4x$.

\end{remark}

\medskip

We would like to cite two such problems in this section: the nature of a very 
particular product on one hand, and a hard general problem on the other hand.

\subsection{The Flajolet-Martin constant}

In a 1985 paper \cite[Theorem~3.1, p.\ 193]{FM} Flajolet and Martin came across the following constant
$$
R := \frac{e^{\gamma}\sqrt{2}}{3}
\prod_{n \geq 1} \left(\frac{(4n+1)(4n+2)}{4n(4n+3)}\right)^{(-1)^{s(n)}}
$$
where $\gamma$ is the Euler constant and, as in the introduction, $s(n)$ is the sum of the binary digits 
of $n$. It is easily proven that $R=\frac{2^{-1/2}e^{\gamma}}{Q}$, where
$$
Q := \prod_{n \geq 1} \left(\frac{2n}{2n+1}\right)^{(-1)^{s(n)}}.
$$
This product resembles the (algebraic) Woods-Robbins product quoted in the introduction very much
and it is not extraordinarily different from the third example (transcendental) given after 
Remark~\ref{k-ell}, respectively
$$
\prod_{n \geq 0} \left(\frac{2n+1}{2n+2}\right)^{(-1)^{s(n)}} = \frac{\sqrt{2}}{2} 
\ \ \ \ \ \mbox{\rm and} \ \ \ \ \
\prod_{n \geq 0} \left(\frac{2n+6}{2n+3}\right)^{\varepsilon_n} =
\frac{\Gamma(1/4)^6}{8\sqrt{2}\,\pi^{7/2}}\cdot
$$
But not only we do not know any closed-form expression for the infinite product $Q$, but also we 
do not know its arithmetical nature, nor the arithmetical nature of the Flajolet-Martin constant.

\begin{remark}\label{do-not-know}
In the same spirit, looking at the infinite products $A$ and $B$ given in the introduction
$$
A = \prod_{n \geq 0} \left(\frac{2n+1}{2n+2}\right)^{\varepsilon_n}
\ \ \ \ \ \mbox{\rm and} \ \ \ \ \
B = \prod_{n \geq 1} \left(\frac{2n}{2n+1}\right)^{\varepsilon_n} = \frac{\Gamma(1/4)^2}{8\sqrt{2 \pi}}
$$
we do not know any closed-form formula for $A$, nor whether this is a transcendental or algebraic number.
\end{remark}

\subsection{The Rohrlich and Rohrlich-Lang conjectures}

A strong conjecture known as the Rohrlich conjecture predicts that the algebraicity of any products and 
quotients of normalized gamma values (the normalized gamma function is $\Gamma/\sqrt{\pi}$) must derive from 
the properties $\Gamma(x+1)=\Gamma(x)$, the reflection formula and the multiplication formula for $\Gamma$ 
(see, e.g., \cite[p.\ 418]{Lang} or \cite[p.\ 444--445]{Waldschmidt}). An explicit formulation is given, e.g., 
in \cite{RR}:

\begin{conjecture}[{\rm Rohrlich}]
Let $a_1, a_2, \ldots, a_r$ be rational numbers that are not in $\{0, -1, -2, \ldots\}$. Let $D$ be a
common denominator of the $a_i$'s. Then the product $\Gamma(a_1) \Gamma(a_2) \cdots \Gamma(a_r)$ is an 
algebraic multiple of $\pi^{r/2}$ if and only if for all $m \in \{1, 2, \ldots, D-1\}$ relatively prime to $D$
we have $\sum_{1 \leq i \leq r} \{ma_i\} = r/2$, where $\{x\}$ denotes the fractional part of the real $x$.
\end{conjecture}

\begin{remark}
There is an even stronger conjecture, known as the Rohrlich-Lang conjecture. The interested reader can
look, e.g., at \cite[p.\ 445]{Waldschmidt}.
\end{remark}

The Rohrlich conjecture implies a criterion for the algebraicity of the infinite products in 
Theorem~\ref{gamma}. We first need an easy lemma.

\begin{lemma}
Let $\{x\}$ denote the fractional part of the real number $x$. Then we have 
$$
\begin{array}{lll}\label{fractional}
\{-x\} &=& 1 - \{x\} \ \mbox{if $x$ is not an integer} \\
\{\frac{1}{2} + x\} &=& 
\begin{cases}
\{x\} + \frac{1}{2} \ &\mbox{\rm if} \ \{x\} < \frac{1}{2} \\
\{x\} - \frac{1}{2} \ &\mbox{\rm if} \ \{x\} \geq \frac{1}{2}
\end{cases} \\
\{\frac{1}{2} - x\} &=&                               
\begin{cases}
\frac{3}{2} - \{x\} \ &\mbox{\rm if} \ \{x\} > \frac{1}{2} \\
\frac{1}{2} - \{x\} \ &\mbox{\rm if} \ \{x\} \leq \frac{1}{2} 
\end{cases}
\end{array}
$$
\end{lemma}

\proof Left to the reader. \endpf

\begin{corollary}

\ { }

\begin{itemize}

\item[(i)] Let $b$ be a positive rational number with denominator $a$ such that $b/4$ is not
an integer. Then, under the conjecture of Rohrlich, we have that the infinite product
$$
\prod_{n \geq 0} \left( \frac{n+b}{n + \frac{1+b}{2}} \right)^{\varepsilon_n}
$$
is algebraic if and only if for all $m \in \{1, 2, \ldots, 4a-1\}$ such that
$m$ is relatively prime to $4a$ we have either
($\{mb/4\} < 1/2$ and $\{m/4\} \leq 1/2$) or
($\{mb/4\} \geq 1/2$ and $\{m/4\} > 1/2$). In other words this condition is equivalent to
saying that for all $m \equiv 1 \bmod 4$ and relatively prime to $a$ we have $\{mb/4\} < 1/2$,
and for all $m \equiv 3\bmod 4$ and relatively prime to $a$ we have $\{mb/4\} \geq 1/2$.

\item[(ii)] Let $b$ and $c$ be positive rational numbers with common denominator $a$. 
We also suppose that neither $b/4$ nor $(c-2)/4$ are integers. Then, under 
the conjecture of Rohrlich, we have that the infinite product
$$
A(b,c) := \prod_{n \geq 0} \left(
\frac{(n+b)(n + \frac{1+c}{2})}{(n+c)(n + \frac{1+b}{2})}
                 \right)^{\varepsilon_n}
$$
is algebraic if and only if for all $m \in \{1, 2, \ldots, 4a-1\}$ such that
$m$ is relatively prime to $4a$ we have either
($\{mb/4\} < 1/2$ and $\{mc/4\} \leq 1/2$) or
($\{mb/4\} \geq 1/2$ and $\{mc/4\} > 1/2$).

\end{itemize}

\end{corollary}

\proof

Since Assertion~(i) is the particular case $c=1$ of Assertion~(ii), it suffices to prove (ii).
Applying Theorem~\ref{gamma} and the reflection formula for the gamma function, we have
$$
A(b,c) = \frac{\Gamma(\frac{c}{4})\Gamma(\frac{1}{2} + \frac{b}{4})}
       {\Gamma(\frac{b}{4})\Gamma(\frac{1}{2} + \frac{c}{4})}
= \frac{1}{\pi^2} \sin \frac{b \pi}{4} \cos \frac{c \pi}{4} \
       \Gamma\left(\frac{c}{4}\right) \Gamma\left(\frac{1}{2} - \frac{c}{4}\right)
       \Gamma\left(\frac{1}{2} + \frac{b}{4}\right) \Gamma\left(1 - \frac{b}{4}\right).   
$$
Since $b$ and $c$ are rational, the product $\sin \frac{b \pi}{4} \cos \frac{c \pi}{4}$ is
algebraic. Hence $A(b,c)$ is algebraic if and only if the quantity
$$
\frac{1}{\pi^2}  \Gamma\left(\frac{c}{4}\right) \Gamma\left(\frac{1}{2} - \frac{c}{4}\right)
       \Gamma\left(\frac{1}{2} + \frac{b}{4}\right) \Gamma\left(1 - \frac{b}{4}\right)
$$
is algebraic (recall the conditions on $b$ and $c$).
Under the conjecture of Rohrlich, we thus see that $A(b,c)$ is algebraic if and only if
for all $m \in \{1, 2, \ldots 4a-1\}$ which is relatively prime to $4a$ we have
$$
\left\{\frac{mc}{4}\right\} + \left\{m\left(\frac{1}{2} - \frac{c}{4}\right)\right\}
+ \left\{m\left(\frac{1}{2} + \frac{b}{4}\right)\right\} + \left\{m\left(1 - \frac{b}{4}\right)\right\}
= 2.
$$
Since $m$ is relatively prime to $4a$, $m$ must be odd, thus the condition becomes
$$
\left\{\frac{mc}{4}\right\} + \left\{\frac{1}{2} - \frac{mc}{4}\right\}
+ \left\{\frac{1}{2} + \frac{mb}{4}\right\} + \left\{-\frac{mb}{4}\right\}
= 2.
$$
Applying Lemma~\ref{fractional} this is equivalent to
$$
\left(\left\{\frac{mb}{4}\right\} < \frac{1}{2} 
\ \mbox{\rm and} \ 
\left\{\frac{mc}{4}\right\} \leq \frac{1}{2}\right)
\ \ \ \mbox{\rm or} \ \ \
\left(\left\{\frac{mb}{4}\right\} \geq \frac{1}{2} 
\ \mbox{\rm and} \ 
\left\{\frac{mc}{4}\right\} > \frac{1}{2}\right)\cdot \ \Box
$$

\section{Conclusion}

The quest for closed-form values is teasing but endless. Since we focussed on the (regular) 
paperfolding sequence here, we would like to cite here a result (due to von Haeseler in the 
case $s=1$, see \cite[exercise~27, p.\ 205]{AS}):
$$
\sum_{n \geq 0} \frac{\varepsilon_n}{(n+1)^s} 
= \frac{2^s}{2^s-1} \sum_{n \geq 0} \frac{(-1)^n}{(2n+1)^s}
\ \ \ \mbox{\rm in particular} \ \ \
\sum_{n \geq 0} \frac{\varepsilon_n}{n+1} = \frac{\pi}{2}\cdot
$$
Note that this result can also be obtained from an additive version of Lemma~\ref{general}.

% http://www.ybook.co.jp/online-p/nao2012pdf/nao2012p239.pdf

\end{document}